\documentclass[10pt]{amsart}
\usepackage{enumerate,amsmath,amssymb,latexsym,
amsfonts, amsthm, amscd, mathabx}

\theoremstyle{plain}

\newtheorem{theorem}{Theorem}

\numberwithin{equation}{section}

\addtocounter{section}{-1}

\begin{document}

\title {Self-dual Harmonicity: the planar case}

\date{}

\author[P.L. Robinson]{P.L. Robinson}

\address{Department of Mathematics \\ University of Florida \\ Gainesville FL 32611  USA }

\email[]{paulr@ufl.edu}

 \keywords{}

\begin{abstract}

We present a manifestly geometrically self-dual version of the Fano harmonicity axiom for the projective plane. 

\end{abstract}

\maketitle

\section*{Introduction} 

In a recent sequence of papers, we have explored manifestly self-dual formulations of three-dimensional projective space: [2] begins the sequence, basing projective space on lines and their abstract incidence, from which points and planes are derived as secondary elements; [4] pursues the complementary approach, taking points and planes to be primary elements from which lines are derived. The most recent paper [6] addresses the axioms (or assumptions) of harmonicity and projectvity: traditionally, {\it harmonicity} asserts (with Fano) that `the diagonal points of a complete quadrangle are noncollinear' while {\it projectivity} asserts that `if a projectivity leaves each of three distinct points of a line invariant, it leaves every point of the line invariant'; see [7] Section 18 and [7] Section 35 respectively. 

\medbreak 

In keeping with our aim to present projective space in a formulation that is manifestly self-dual, it is desirable to express projectivity and harmonicity in manifestly self-dual forms. The axiom of projectivity is so expressed in [3] and [6]: indeed, [7] Section 103 contains a self-dual expression of projectivity in its discussion of reguli. In [5] and [6] our versions of the axiom of harmonicity are self-dual only in being the conjunction of a pair of dual statements, each of which implies the other; we might say that these versions are self-dual in rather a logical sense. Our purpose in the present paper and its sequel is to offer expressions of harmonicity that are more geometrically self-dual: the present paper deals with harmonicity in the plane; its sequel deals with the more elaborate spatial harmonicity. 

\medbreak 

\section*{The Planar Axiom}

\medbreak 

In this section, we present a manifestly self-dual version of the harmonicity axiom in the projective plane. Recall that {\it harmonicity} is traditionally formulated as follows: 
\medbreak 
`the diagonal points of a complete quadrangle are noncollinear'
\medbreak 
\noindent
and that the planar dual of this reads as follows: 
\medbreak 
`the diagonal lines of a complete quadrilateral are noncurrent'.
\medbreak 
\noindent
These two dual statements are equivalent within the Veblen-Young axiomatization [7]: granted the Veblen-Young assumptions of alignment and extension, each statement  implies the other. Of course, we may take the conjunction of this dual pair of equivalent statements as a `self-dual' formulation of the planar harmonicity axiom. It might be objected that the self-dual nature of this version is as much logical as it is geometrical. Our purpose here is to offer a version in which the self-duality is more purely geometrical. In order to present this alternative version, we recall some elementary facts pertaining to quadrangles and quadrilaterals.  

\medbreak

Explicitly, let $P_0 P_1 P_2 P_3$ be a {\it complete quadrangle} in the plane: thus, its four vertices $P_0,  P_1,  P_2,  P_3$ are in general position, in the sense that no three of them are collinear; by definition, the three {\it diagonal points} of this quadrangle are 
$$D_1 = P_0 P_1 \cdot P_2 P_3, \; D_2 = P_0 P_2 \cdot P_3 P_1, \; D_3 = P_0 P_3 \cdot P_1 P_2.$$

\medbreak 

Dually, let $p_0 p_1 p_2 p_3$ be a {\it complete quadrilateral} in the plane: its four sides $p_0, p_1, p_2, p_3$ are in general position, in the sense that no three of them are concurrent; by definition, the three {\it diagonal lines} of this quadrilateral are 
$$d_1 = (p_0 \cdot p_1)(p_2 \cdot p_3), \; d_2 = (p_0 \cdot p_2)( p_3 \cdot p_1), \; d_3 = (p_0 \cdot p_3)(p_1 \cdot p_2).$$

\medbreak 

We note that the diagonal points of a complete quadrangle are distinct: indeed, were we to suppose that $D_1 = D_2 = D$ say, then $D$ would lie on $P_0 P_1, P_0 P_2, P_2 P_3, P_3 P_1$ and hence support the absurd conclusion 
$$P_0 = P_0 P_1 \cdot P_0 P_2 = D = P_2 P_3 \cdot P_3 P_1 = P_3.$$ Dually, the diagonal lines of a complete quadrilateral are distinct. 
\medbreak

Now, let let us start with a complete quadrangle $P_0 P_1 P_2 P_3$. We shall {\it assume} that its diagonal points $D_1, D_2, D_3$ are noncollinear; they are therefore the vertices of its {\it diagonal triangle} $D_1 D_2 D_3$. By the very definition of the diagonal points, the triangles $D_1 D_2 D_3$ and $P_1 P_2 P_3$ are perspective from $P_0$: by the Desargues theorem, they are therefore perspective from a line $p_0$; thus, 
$$p_0 \; \ni \; D_2 D_3 \cdot P_2 P_3, \; D_3 D_1 \cdot P_3 P_1, \; D_1 D_2 \cdot P_1 P_2. $$
The diagonal triangle $D_1 D_2 D_3$ is similarly perspective to $P_0 P_3 P_2$ (from $P_1$), $P_3 P_0 P_1$ (from $P_2$), $P_2 P_1 P_0$ (from $P_3$); accordingly, the Desargues theorem likewise furnishes respective lines
$$p_1 \; \ni \; D_2 D_3 \cdot P_3 P_2, \; D_3 D_1 \cdot P_2 P_0, \; D_1 D_2 \cdot P_0 P_3,$$
$$p_2 \; \ni \; D_2 D_3 \cdot P_0 P_1, \; D_3 D_1 \cdot P_1 P_3, \; D_1 D_2 \cdot P_3 P_0,$$
$$p_3 \; \ni \; D_2 D_3 \cdot P_1 P_0, \; D_3 D_1 \cdot P_0 P_2, \; D_1 D_2 \cdot P_2 P_1.$$
\medbreak
\noindent
The four lines $p_0,  p_1,  p_2,  p_3$ so constructed constitute a complete quadrilateral. For example, suppose $p_1, p_2, p_3$ were to concur at $Q$: then $Q = p_2 \cdot p_3 = D_2 D_3 \cdot P_0 P_1$ would lie on $D_2 D_3$ and $Q = p_3 \cdot p_1 = D_3 D_1 \cdot P_0 P_2$ would lie on $D_3 D_1$ so that $Q = D_3$ by noncollinearity of the diagonal points; but then $D_3 = Q = p_1 \cdot p_2 = D_1 D_2 \cdot P_0 P_3$ would lie on $D_1 D_2$ and noncollinearity of the diagonal points would be contradicted. From above, $p_0 \cdot p_3 = D_1 D_2 \cdot P_1 P_2$ and $p_1 \cdot p_2 = D_1 D_2 \cdot P_0 P_3$: thus $D_1 D_2$ contains both $p_0 \cdot p_3$ and $p_1 \cdot p_2$ and so coincides with the diagonal line $d_3 = (p_0 \cdot p_3)(p_1 \cdot p_2)$ of the quadrilateral; likewise, $D_2 D_3 = d_1$ and $D_3 D_1 = d_2$. 

\medbreak 

Incidentally, we feel free here to use the Desargues theorem in the plane, because our ultimate concern is with projective space (in which all planes are automatically Desarguesian); of course, we might instead have adopted a framework for the projective plane that axiomatically incorporates the Desargues theorem, as in [1]. 

\medbreak 

Thus, each complete quadrangle $P_0 P_1 P_2 P_3$ with noncollinear diagonal points engenders a complete quadrilateral $p_0 p_1 p_2 p_3$ sharing the very same diagonal triangle. Dually, each complete quadrilateral with nonconcurrent diagonal lines engenders a complete quadrangle sharing the very same diagonal triangle. In fact, these dual constructions are mutually inverse; application of the second construction to the quadrilateral that results from the first construction reproduces the original quadrangle $P_0 P_1 P_2 P_3$. To see this, by symmetry we need only show that if $\{ i, j, k \} = \{ 1, 2, 3 \}$ then the line joining $d_i \cdot d_j$ (a vertex of the triangle $d_1 d_2 d_3$) and $p_j \cdot p_j$ (the corresponding vertex of $p_1 p_2 p_3$) passes through $P_0$; but this is clear, for $p_j \cdot p_j = D_i D_j \cdot P_k P_0$ lies on $P_k P_0$ and $d_i \cdot d_j = D_j D_k \cdot D_k D_i = D_k \notin D_i D_j$ lies on $P_k P_0$ so that $(d_i \cdot d_j)(p_j \cdot p_j) = P_k P_0$. 

\medbreak 

The figure comprising a complete quadrangle and a complete quadrilateral that share the same diagonal triangle is called a {\it quadrangle-quadrilateral configuration}. We may present this configuration in a way that does not make quite explicit the fact that the respective diagonals form a triangle. Thus, the configuration comprises a complete quadrangle $P_0 P_1 P_2 P_3$ (with distinct diagonal points $D_1, D_2, D_3$) and a complete quadrilateral $p_0 p_1 p_2 p_3$ (with distinct diagonal lines $d_1, d_2, d_3$) so related that 
$$d_1 = D_2 D_3, \; d_2 = D_3 D_1, \; d_3 = D_1 D_2$$
and 
$$D_1 = d_2 \cdot d_3, \; D_2 = d_3 \cdot d_1, \; D_3 = d_1 \cdot d_2.$$
\medbreak 
\noindent
Notice that these two sets of relations are equivalent: if $\{ i, j, k \} = \{ 1, 2, 3 \}$ and the former set holds, then the distinct lines $d_i$ and $d_j$ have $D_k$ in common, so $D_k = d_i \cdot d_j$ and the latter set holds. As we note in the following paragraph, these relations imply that $D_1, D_2, D_3$ are noncollinear and that $d_1, d_2, d_3$ are nonconcurrent, so that $P_0 P_1 P_2 P_3$ and $p_0 p_1 p_2 p_3$ indeed share the same diagonal triangle. For more on the quadrangle-quadrilateral configuration, see the account in [7] Section 18; for a slightly different perspective on the construction that leads from a complete quadrangle to the complete quadrilateral sharing its diagonal triangle, see also [1] Section 3.1 Exercise 2. 

\medbreak 

To summarize the foregoing discussion, if the diagonal points of a complete quadrangle $P_0 P_1 P_2 P_3$ are not collinear, then there exists a corresponding complete quadrilateral $p_0 p_1 p_2 p_3$ making up a quadrangle-quadrilateral configuration; dually,  if the diagonal lines of a complete quadrilateral $p_0 p_1 p_2 p_3$ are not concurrent, then there exists a corresponding complete quadrangle $P_0 P_1 P_2 P_3$ making up a quadrangle-quadrilateral configuration. In the opposite direction, if a complete quadrangle $P_0 P_1 P_2 P_3$ (with distinct diagonal points $D_1, D_2, D_3$) and a complete quadrilateral $p_0 p_1 p_2 p_3$ (with distinct diagonal lines $d_1, d_2, d_3$) are so related as to satisfy the relations 
$$d_1 = D_2 D_3, \; d_2 = D_3 D_1, \; d_3 = D_1 D_2$$
and 
$$D_1 = d_2 \cdot d_3, \; D_2 = d_3 \cdot d_1, \; D_3 = d_1 \cdot d_2$$
\medbreak 
\noindent
then the diagonal points $D_1, D_2, D_3$ of $P_0 P_1 P_2 P_3$ are not collinear and the diagonal lines $d_1, d_2, d_3$ of $p_0 p_1 p_2 p_3$ are not concurrent: for instance, $D_1 \neq D_2$ implies that $d_2 \cdot d_3 \neq d_3 \cdot d_1$ so that $d_1, d_2, d_3$ do not concur. 

\medbreak 

Recalling that if one complete quadrangle has noncollinear diagonal points then so do all, we have now established the following alternative self-dual version of the planar harmonicity axiom: 
\medbreak 
[{\bf H}] There exists a quadrangle-quadrilateral configuration. 

\medbreak 

Notice that the self-dual nature of this version is manifestly geometrical. Of course, this version is equivalent both to the assertion that each quadrangle/quadrilateral is part of a quadrangle-quadrilateral configuration and to the original version according to which the diagonals of each quadrangle/quadrilateral form a triangle. 

\medbreak 

It is of some interest to observe that we may interrupt the present discussion at an intermediate point. Again let $P_0 P_1 P_2 P_3$ be a complete quadrangle with noncollinear diagonal points $D_1, D_2, D_3$: our construction furnishes a complete quadrilateral $p_0 p_1 p_2 p_3$ such that if $\{ i, j, k \} = \{ 1, 2, 3 \}$ then $p_i \cdot p_j = P_0 P_k \cdot D_i D_j$ and $p_0 \cdot p_k = P_i P_j \cdot D_i D_j$; of course, the sides $p_0, p_1, p_2, p_3$ themselves can be recovered from these points. Our dual construction leads from a complete quadrilateral $p_0 p_1 p_2 p_3$ with nonconcurrent diagonal lines $d_1, d_2, d_3$ to a complete quadrangle $P_0 P_1 P_2 P_3$ whose vertices can be recovered from the lines $P_i P_j = (p_0 \cdot p_k)(d_i \cdot d_j)$ and $P_0 P_k = (p_i \cdot p_j)(d_i \cdot d_j)$. Temporarily suspending the harmonicity axiom, we may say that the complete quadrangle $P_0 P_1 P_2 P_3$ (with diagonal points $D_1, D_2, D_3$) and the complete quadrilateral $p_0 p_1 p_2 p_3$ (with diagonal lines $d_1, d_2, d_3$) are {\it mated} iff they satisfy the relations 
$$p_i \cdot p_j = P_0 P_k \cdot D_i D_j, \; p_0 \cdot p_k = P_i P_j \cdot D_i D_j$$
and 
$$P_i P_j = (p_0 \cdot p_k)(d_i \cdot d_j), \; P_0 P_k = (p_i \cdot p_j)(d_i \cdot d_j)$$
\medbreak 
\noindent 
whenever $\{ i, j, k \} = \{ 1, 2, 3 \}$. It may be readily verified that these two sets of conditions are both equivalent to each other and equivalent to the condition that $P_0 P_1 P_2 P_3$ and $p_0 p_1 p_2 p_3$ make up a quadrangle-quadrilateral configuration. 

\medbreak 

Accordingly, the following is another manifestly self-dual version of planar harmonicity: 
\medbreak 
[${\bf H}'$] There exists a mated quadrangle-quadrilateral pair. 

\medbreak 

Equivalently, each quadrangle/quadrilateral has a mate. 

\medbreak

\bigbreak 

\begin{center} 
{\small R}{\footnotesize EFERENCES}
\end{center} 
\medbreak 

[1] H.S.M. Coxeter, {\it Projective Geometry}, Second Edition, Springer-Verlag (1987). 

\medbreak 

[2] P.L. Robinson, {\it Projective Space: Lines and Duality}, arXiv 1506.06051 (2015). 

\medbreak 

[3] P.L. Robinson, {\it Projective Space: Reguli and Projectivity}, arXiv 1506.08217 (2015). 

\medbreak 

[4] P.L. Robinson, {\it Projective Space: Points and Planes}, arXiv 1611.06852 (2016). 

\medbreak 

[5] P.L. Robinson, {\it Projective Space: Tetrads and Harmonicity}, arXiv 1612.01913 (2016). 

\medbreak 

[6] P.L. Robinson, {\it Projective Space: Harmonicity and Projectivity}, arXiv 1612.08422 (2016). 

\medbreak 

[7] O. Veblen and J.W. Young, {\it Projective Geometry}, Volume I, Ginn and Company (1910).

\medbreak

\end{document}